\documentclass[11pt]{article}
\usepackage{a4wide}
\usepackage{amsmath}
\usepackage{amsfonts}
\usepackage{iwona}
\usepackage{latexsym}
\usepackage{amssymb}

\date{}

\parskip\medskipamount
\newtheorem{Theorem}{Theorem}

\newtheorem{Lemma}{Lemma}

\title{A Monomial Basis for the Holomorphic Functions on $c_{0}$}

\author{Se\'an Dineen and Jorge Mujica
\footnote{MSC 2010 Classification 46G20, 32A05.
 }}
\begin{document}
\maketitle 
\allowdisplaybreaks

\begin{abstract}
  We show that both $(\mathcal{H}(c_{0}),\tau_{\omega})$ and
  $(\mathcal{H}_{b}(c_{0}),\tau_{b})$ have a monomial Schauder basis.
\end{abstract}

\section*{Introduction} Various authors have considered
monomial expansions of polynomials defined on infinite dimensional
Banach spaces (\cite{Alencar, Defant, Dimant, Grecu, Matos, Ryan})
and although it has been shown that they form a Schauder basis in
some cases so far we do not have analogous results for spaces of
holomorphic functions. To show that the monomials form a basis for
various spaces of holomorphic functions on $c_{0}$ we use three
different decompositions; $\mathcal{S}_{*}$-absolute
decompositions of locally convex spaces \cite{Dineen1}, finite
dimensional monotone decompositions of a Banach space
\cite{Dimant}, and a Schauder basis \cite{Diestel}. We discuss
these in section~1. In section~2 we recall the definitions of the
different spaces of holomorphic functions and discuss the square
order on the polynomials and prove our main result in section~3.
\section{Linear Decompositions} A sequence of subspaces $\{E_{n}\}_{n=1}^{\infty}$ of a
locally convex space $E$ is called a decomposition for $E$ if for
each $x\in E$ there exists a unique sequence
$(x_{n})_{n=1}^{\infty}$,$x_{n}\in E_{n}$ for all $n$ such that
\begin{eqnarray*}x=\sum_{n=1}^{\infty}x_{n}:=\lim_{n\to
\infty}\sum_{j=1}^{n}x_{j}.\end{eqnarray*} If $E$ admits a
fundamental system of semi-norms $\mathcal{N}$ such that for any
$p\in \mathcal{N}$ and any sequence of scalars
$(\alpha_{n})_{n=1}^{\infty}$ satisfying $\limsup_{n\to
\infty}|\alpha_{n}|^{1/n}<\infty$ the semi-norm $q$ defined as
follows:
\begin{eqnarray}q\big{(}\sum_{n=1}^{\infty}x_{n}\big{)}:=\sum_{n=1}^{\infty}|\alpha_{n}|p(x_{n})
\end{eqnarray} is continuous, then we say that $\{E_{n}\}_{n=1}^{\infty}$ is an $\mathcal{S}_{*}$-absolute
decomposition of $E$. This concept coincides with the notion of
global Schauder decomposition given in \cite{Galindo} and
\cite{Venkova} and is a variation on $\mathcal{S}$-absolute
decomposition discussed in \cite{Dineen1}.

If each $E_{n}$ is a finite dimensional space the decomposition is
called finite dimensional. If $E$ is a normed linear space with
norm $\|\cdot\|$ then the infimum of all $c$ such that
\begin{eqnarray}\|\sum_{j=1}^{n}x_{j}\|\leq c\|\sum_{j=1}^{m}x_{j}\|\end{eqnarray} for all
positive integers $m$ and $n$, $n<m$, and all $x_{j}\in E_{j}$ is
called the decomposition constant. The decomposition constant is
always greater than or equal to 1 and if it equals 1 we say that
the decomposition is monotone. A renorming generally changes the
decomposition constant.

If each $E_{n}$ is one dimensional and $e_{n}$ spans $E_{n}$ we
say that $(e_{n})_{n=1}^{\infty}$ is a Schauder basis for $E$. In
this case there exists for each $x$ a sequence of scalars
$(x_{n})_{n=1}^{\infty}$ such that
$x=\sum_{n=1}^{\infty}x_{n}e_{n}$ and we use the term basis
constant in place of decomposition constant. 
The linear functional $x\longrightarrow x_{n}$ is called the
$n^{th}$ coefficient functional and is denoted by $e_{n}^{*}$.

If $E$ has an $\mathcal{S}_{*}$-decomposition,
$\{E_{n}\}_{n=1}^{\infty}$, then it admits, by (1), a fundamental
system of semi-norms $\mathcal{N}$ such that
\begin{eqnarray}p\big{(}\sum_{n=1}^{\infty}x_{n}\big{)}=\sum_{n=1}^{\infty}p(x_{n})\end{eqnarray}
for all $\sum_{n=1}^{\infty}x_{n}\in E, x_{n}\in E_{n}$ for all
$n$. If each $E_{n}$ has a Schauder basis
$(e_{n,m})_{m=1}^{\infty}$ then an ordering of
$(e_{n,m})_{n,m=1}^{\infty}$ into a sequence is given by a
bijective mapping $\phi:\mathbb{N}^{2}\longrightarrow\mathbb{N}$.
We say that the ordering is compatible if
\begin{eqnarray}m<\overline{m}\Longrightarrow\phi(n,m)<\phi(n,\overline{m})\end{eqnarray}
for all $n$, that is it induces on each $E_{n}$ its original order
(see Proposition~4.1 in \cite{Dineen1}).
\begin{Lemma} Let $\phi:\mathbb{N}^{2}\longrightarrow\mathbb{N}$ denote a compatible
ordering. For every positive integer $j$, there exists a finite
subset $S_{j}$ of positive integers, and a finite set of positive
integers $(k_{n}(j))_{n\in S_{j}}$ such that
\begin{eqnarray}\{k:1\leq k\leq j\}=\bigcup_{n\in
S_{j}}\{\phi(n,m):1\leq m\leq k_{n}(j)\}.\end{eqnarray} Moreover,
if $l$ is a positive integer then $S_{j+l}=S_{j}\cup S$ for some
finite subset $S$ of $\mathbb{N}$, disjoint from $S_{j}$, and
$k_{n}(j)\leq k_{n}(j+l)$ for all $n\in S_{j}$.

\end{Lemma}{\bf{Proof.}} We prove this result by induction on $j$. Since $\phi$ is surjective
there exists a pair of integers $(n_{1},m_{1})$ such that
$\phi(n_{1},m_{1})=1$. Let $S_{1}=\{n_{1}\}$. If $m_{1}>1$ then,
since $\phi$ is injective we have
$\phi(n_{1},1)>1=\phi(n_{1},m_{1})$ and this contradicts (4).
Hence $m_{1}=1$ and (5) holds when $j=1$.

Now suppose (5) holds for the positive integer $j$. This implies,
in particular, that \begin{eqnarray}\{1,2,\ldots,j\}\subset
\bigcup_{n\in S_{j}}\{\phi(n,m):m \in N \}\textrm{
}.\end{eqnarray} By surjectivity of $\phi$ we can find a pair of
positive integers $(n_{j+1},m_{j+1})$ such that
$\phi(n_{j+1},m_{j+1})=j+1$. We consider two cases.

If $n_{j+1}\not\in S_{j}$ then, since $\phi$ is injective, (6)
implies that $\phi(n_{j+1},m)\geq j+1$ for all $m\in N$. If
$m_{j+1}>1$ then, since $\phi(n_{j+1},m_{j+1})=j+1$, the
injectivity of $\phi$ implies $\phi(n_{j+1},1)>j+1$. We then have
$\phi(n_{j+1},1)>j+1=\phi(n_{j+1},m_{j+1})$ and this contradicts
(4). Hence $m_{j+1}=1$. Letting $S_{j+1}=S_{j}\cup \{n_{j+1}\}$,
$k_{n}(j)=k_{n}(j+1)$ for $n\in S_{j}$ and $k_{n_{j+1}}(j+1)=1$ we
obtain (5).

If $n_{j+1}\in S_{j}$ let $S_{j+1}=S_{j}$. If
$m_{j+1}>k_{n_{j+1}}(j)+1$ then, by (6) and since $\phi$ is
injective,
$$\phi(n_{j+1},k_{n_{j+1}}(j)+1)>j+1=\phi(n_{j+1},m_{j+1})$$ and
this contradicts (4). Hence $m_{j+1}\leq k_{n_{j+1}}(j)+1$. By
(6), $\phi(n_{j+1},l)\leq j$ for all $l\leq k_{n_{j+1}}(j)$ and,
as $\phi$ is injective, this implies $m_{j+1}>k_{n_{j+1}}(j)$.
Hence $m_{j+1}=k_{n_{j+1}}(j)+1.$ If $n\in S_{j}, n\not=n_{j+1},$
let $k_{n}(j+1)=k_{n}(j)$ and let
$k_{n_{j+1}}(j+1)=k_{n_{j+1}}(j)+1.$
 This implies (5) holds for $j+1$. By induction this completes the proof of
(5) and the remainder of the proof follows easily.

\begin{Theorem} Suppose $E$ has an $\mathcal{S}_{*}$-decomposition,
$\{E_{n}\}_{n=1}^{\infty}$, with fundamental system of semi-norms
$\mathcal{N}$ satisfying (3) and that each $E_{n}$ has a Schauder
basis $(e_{n,m})_{m=1}^{\infty}$. For each $p\in \mathcal{N}$ and
each positive integer $n$ let $c_{p,n}$ denote the basis constant
for $(E_{n},p)$. If for each $p\in\mathcal{N}$, $limsup_{n\to
\infty}c_{p,n}^{1/n}<\infty$ then $(e_{n,m})_{n,m=1}^{\infty}$,
with any compatible ordering, is a basis for $E$.
\end{Theorem}
{\bf{Proof.}} Let $\phi:\mathbb{N}^{2}\longrightarrow\mathbb{N}$
denote a fixed compatible order on $\mathbb{N}^{2}$. By the
definition of decomposition and basis we see that
$(e_{n,m})_{n,m=1}^{\infty}$ spans a dense subspace of $E$ and
hence if suffices to show it is a basic sequence in E. To show
this we apply Theorem 6, p. 298, in \cite{Jarchow}. Let $p$ denote
a semi-norm on $E$ satisfying (3). If $c_{p,n}$ is the basis
constant for $(E_{n},p)$ then, by (1), the semi-norm
$$q\big{(}\sum_{k=1}^{\infty}x_{n}\big{)}:=\sum_{k=1}^{\infty}c_{p,n}p\big{(}x_{n}\big{)}, x_{n}\in E_{n},
\sum_{n=1}^{\infty}x_{n}\in E$$ is continuous on $E$.

 Let $(\alpha_{k})_{k\in
\mathbb{N}}$ denote an arbitrary set of scalars. We now used the
notation employed in Lemma~1. If $j$ and $l$ are positive integers
then $S_{j+l}=S_{j}\cup S$ for some finite subset $S\subset
\mathbb{N}$ disjoint from $S_{j}$, and $k_{n}(j)\leq k_{n}(j+l)$
for all $n\in S_{j}$. We then have for all positive integers $j$
and $l$,
\begin{align*}p\big{(}\sum_{k=1}^{j}\alpha_{k}e_{\phi^{-1}(k)}\big{)}
&= p\big{(}\sum_{n\in S_{j}}\big{\{}
\sum_{m=1}^{k_{n}(j)}\alpha_{\phi(n,m)}
e_{n,m}\big{\}}\big{)}\\
&=\sum_{n\in S_{j}}p\big{(}
\sum_{m=1}^{k_{n}(j)}\alpha_{\phi(n,m)} e_{n,m}\big{)}\\
&\leq \sum_{n\in S_{j}}c_{p,n}p\big{(}
\sum_{m=1}^{k_{n}(j+l)}\alpha_{\phi(n,m)} e_{n,m}\big{)}\\
&\leq \sum_{n\in S_{j}}c_{p,n}p\big{(}
\sum_{m=1}^{k_{n}(j+l)}\alpha_{\phi(n,m)}
e_{n,m}\big{)} \\
&\qquad\qquad +c_{n,p}p\big{(}\sum_{n\in S}\big{\{}
\sum_{m=1}^{k_{n}(j+l)}\alpha_{\phi(n,m)}
e_{n,m}\big{\}}\big{)}\\
&=\sum_{n\in S_{j+l}}c_{p,n}p\big{(}
\sum_{m=1}^{k_{n}(j+l)}\alpha_{\phi(n,m)}
e_{n,m}\big{)}\\
&= q\big{(}\sum_{k=1}^{j+l}\alpha_{k}e_{\phi^{-1}(k)}\big{)}\end{align*}
This completes the proof.

 Let $c_{0}=\{(z_{j})_{j=1}^{\infty}\subset \mathbb{C}: \lim_{j\to \infty} z_{j}=0\}$
and let $c_{0}^{+}=\{(z_{j})_{j=1}^{\infty}\in c_{0}:z_{j}>0\}$.
We denote by $(e_{j})_{j=1}^{\infty}$ the standard unit vector
basis for $c_{0}$ and let $(e_{j}^{*})_{j=1}^{\infty}$ denote the
dual unit vector basis for $\ell_{1}=c_{0}^{\prime}.$ The
polydiscs
$$P_{z}:=P_{(z_{j})_{j=1}^{\infty}}=\{(w_{j})_{j=1}^{\infty}\in
c_{0}:|w_{j}|\leq z_{j}\textrm { all } j\}$$ form a fundamental
system of compact subsets for $c_{0}$ when
$z=(z_{j})_{j=1}^{\infty}$ ranges over $c_{0}^{+}$ while the
polydiscs $(B_{r}:=\{z\in c_{0}:\|z\|< r\})_{r=1}^{\infty}$ are a
fundamental system of bounded subsets of $c_{0}$.

\section{Polynomials and Holomorphic Functions} 
In this section we discuss concepts from
infinite dimensional holomorphy and refer to \cite{Dineen1} and
\cite{Mujica} for details. Our main result concerns holomorphic
functions on $c_{0}$.

For each positive integer $n$ let $\mathcal{P}(^{n}c_{0})$ denote
the space of continuous $n$-homogeneous polynomials on $c_{0}$.
Endowed with the supremum norm of uniform convergence over the
unit ball of $c_{0}$, $\mathcal{P}(^{n}c_{0})$ is a Banach space.
Let $\mathbb{N}^{(\mathbb{N})}$ denote the set of all sequences of
non-negative integers which are eventually zero. If
$(m_{i})_{i=1}^{\infty}\in \mathbb{N}^{(\mathbb{N})}$ we call
$|m|:=\sum_{i}m_{i}$ and $l(m):=\sup\{i: m_{i}\not=0\}$ the
modulus and length of $m$, respectively, and call the mapping
$$(z_{j})_{j=1}^{\infty}\in c_{0}\longrightarrow
z^{m}:=z_{1}^{m_{1}}\cdots z_{s}^{m_{s}}\cdots $$ a monomial (we
use the convention $0^{0}=1$). For positive integers $n$ and $k$
let $\mathcal{P}_{k}(^{n}c_{0})$ denote the subspace of
$\mathcal{P}(^{n}c_{0})$ spanned by $\{z^{m}:l(m)=k,|m|=n\}$. For
all $n$ the sequence
$\{\mathcal{P}_{k}(^{n}c_{0})\}_{k=1}^{\infty}$ is a finite
dimensional decomposition of $\mathcal{P}(^{n}c_{0})$ (see
\cite{Dimant} and \cite{Dineen1}, section 4.1). This follows, for
$P\in \mathcal{P}(^{n+1}c_{0})$, from the identities
\begin{eqnarray*}P(\sum_{j=1}^{\infty}z_{j}e_{j})
&=&P(z_{1}e_{1})+\sum_{k=1}^{\infty}\big{\{}P(\sum_{j=1}^{k+1}z_{j}e_{j})-P(\sum_{j=1}^{k}z_{j}e_{j})\big{\}}
\\&=&a_{1}z_{1}^{n+1}+\sum_{s=1}^{n+1} a_{s}z_{1}^{n+1-s}z_{2}^{s}+\sum_{s=1,t\geq 0,s+t<n+1}^{n+1}
 a_{s,t}z_{1}^{n-s-t}z_{2}^{t}z_{3}^{s}+\cdots\end{eqnarray*}
and
\begin{eqnarray*}Q_{k+1}(\sum_{j=1}^{\infty}z_{j}e_{j}):=P(\sum_{j=1}^{k+1}z_{j}e_{j})
 -P(\sum_{j=1}^{k}z_{j}e_{j})=:R_{k+1}(\sum_{j=1}^{\infty}z_{j}e_{j})\cdot
 z_{k+1}\end{eqnarray*} where $Q_{k+1}\in \mathcal{P}_{k+1}(^{n+1}c_{0})$ and
 $R_{k+1}\in \bigoplus_{j=1}^{k+1}\mathcal{P}_{j}(^{n}c_{0}).$
Hence \begin{eqnarray*}Q_{k+1}=R_{k+1}\cdot
e_{k+1}^{*}.\end{eqnarray*}
 We consider different equivalent norms on $P\in
 \mathcal{P}(^{n}c_{0})$, generated by uniform
 convergence over bounded polydiscs in $c_{0}$. If
 $$A:=\{(z_{j})_{j=1}^{\infty}:|z_{j}|\leq \lambda_{j},j=1,2,\ldots\}$$ is a bounded polydisc in
 $c_{0}$, and $Q_{k}\in
 \mathcal{P}_{k}(^{n}c_{0})$ for all $k$, then for all $s,t$ with $s<t,$ we
 have, since $Q_{l}\big{(}\sum_{j=1}^{k}z_{j}e_{j}\big{)}=0$ for
 all $l>k$,
 \begin{eqnarray}\|\sum_{k=1}^{s}Q_{k}\|_{A}\leq
 \|\sum_{k=1}^{t}Q_{k}\|_{A}\end{eqnarray} and this implies that
 $\{\mathcal{P}_{k}(^{n}c_{0})\}_{k=1}^{\infty}$ is a finite dimensional
 monotone decomposition of $\mathcal{P}(^{n}c_{0})$ for all $n$. Using the notation in
 (6) and the fact that holomophic functions on polydiscs achieve their absolute
 maxima on the distinguished boundary we see that
 \begin{eqnarray}\|Q_{k+1}\|_{A}=\lambda_{k+1}\|R_{k+1}\|_{A}.\end{eqnarray}

 We now define the square order on
 $\mathcal{P}(^{n}c_{0})$. On $c_{0}^{\prime}=\ell_{1}$ we use the sequential order inherited from the
 standard unit vector basis $(e_{j}^{*})_{j=1}^{\infty}$. The square order on
 $\mathcal{P}(^{n+1}c_{0})$ is defined as
 follow:
 if $m=(m_{i})_{i=1}^{\infty}$ and $m^{\prime}=(m^{\prime}_{i})_{i=1}^{\infty}$ are in $\mathbb{N}^{(\mathbb{N})}$
 and $|m|=|m^{\prime}|$ then $m<m^{\prime}$ if either $l(m)<l(m^{\prime})$
 or $l(m)=l(m^{\prime})$ and for some positive integer $s\leq l(m)$, $m_{s}<m^{\prime}_{s}$ and $m_{t}=m^{\prime}_{t}$ for
 all $t>s$.

 The square order appears naturally when we use the finite dimensional decomposition
 $(\mathcal{P}_{k}(^{n}c_{0}))_{k=0}^{\infty}$. Clearly, if $k<k^{\prime}$
 then the monomials in $\mathcal{P}_{k}(^{n+1}c_{0})$ precede those
 in $\mathcal{P}_{k^{\prime}}(^{n+1}c_{0})$ and the order within $\mathcal{P}_{k}(^{n+1}c_{0})$ is
 determined by the order inherited from $\mathcal{P}(^{n}c_{0})$. If $m \in \mathbb{N}^{(\mathbb{N})}$ and $l(m)=s$ then
 there is a unique $\underline{m}\in \mathbb{N}^{(\mathbb{N})}$
 such that $z^{m}=z^{\underline{m}}z_{s}$ for all $\sum_{j=1}^{\infty}z_{j}e_{j}\in c_{0}$.
 Note that $|\underline{m}|=|m|-1$ and $l(\underline{m})\leq l(m)$. If $m,m^{\prime}\in
 \mathbb{N}^{(\mathbb{N})}$ then $m<m^{\prime}$ if either
 $l(m)<l(m^{\prime})$ or $l(m)=l(m^{\prime})$ and
 $\underline{m}<\underline{m^{\prime}}$.

The square order was introduced by Ryan (\cite{Ryan}) and various
authors have shown that the monomials of degree $n$ with the
square order are a Schauder basis for $\mathcal{P}(^{n}c_{0})$.
Theorem 2 contains within it yet another proof of this fact,
modulo the result of W. Bogdanowicz and A.Pelczy\'{n}ski in 1957
(see \cite{Dineen1} p. 81) that polynomials on $c_{0}$ are weakly
continuous on bounded sets. The order is important as Defant and
Kalton have shown in \cite{Defant} that when the monomials of
degree $n$ for any $ n\geq 2$ form an unconditional basis for the
space of $n$-homogeneous polynomials on the Banach space $X$,
endowed with the norm of uniform convergence over the unit ball of
$X$, then $X$ is finite dimensional.

We let $\mathcal{H}(c_{0})$ denote the space of holomorphic
functions on $c_{0}$ endowed with the ported $\tau_{\omega}$
topology of Nachbin and let $\mathcal{H}_{b}(c_{0})$ denote the
subspace of $\mathcal{H}(c_{0})$ consisting of all $f$ bounded on
bounded subsets of $c_{0}$ endowed with the Fr\'{e}chet locally
convex topology, $\tau_{b}$, generated by uniform convergence over
the bounded subsets of $c_{0}$. A fundamental system of semi-norms
for $(\mathcal{H}(c_{0}),\tau_{\omega})$ is given by:
\begin{eqnarray}p(\sum_{n=0}^{\infty}P_{n}):=\sum_{n=0}^{\infty}\sup\{|P_{n}(z)|:z=(z_{j})_{j=1}^{\infty}\in
c_{0}, |z_{j}|\leq \beta_{j}+\alpha_{n}\}\end{eqnarray} where
$\sum_{n=0}^{\infty}P_{n}\in \mathcal{H}(c_{0}), P_{n}\in
\mathcal{P}(^{n}c_{0})$ all $n$, and $(\beta_{j})_{j=1}^{\infty}$
and $(\alpha_{n})_{n=0}^{\infty}$ range over $c_{0}^{+}$. Using
the Taylor series expansion at the origin we see that
$(\mathcal{P}(^{n}c_{0}))_{n=0}^{\infty}$ is an
$\mathcal{S}_{*}$-absolute decomposition for both
 $(\mathcal{H}_{b}(c_{0}),\tau_{b})$ and $(\mathcal{H}(c_{0}),\tau_{\omega}).$

 \section{A Schauder Basis for $(\mathcal{H}(c_{0}),\tau_{\omega})$ and $(\mathcal{H}_{b}(c_{0}),\tau_{b})$}
In this section we let $(P_{n,m})_{m=1}^{\infty}$ denote the
monomials of degree $n$ on $c_{0}$ endowed with the square order
and we suppose that the set of all monomials are given a
compatible order.
\begin{Theorem} The monomials with a compatible order are a
Schauder basis for
 $(\mathcal{H}(c_{0}),\tau_{\omega})$ and $(\mathcal{H}_{b}(c_{0}),\tau_{b})$.\end{Theorem}{\bf{Proof.}}
 Both proofs are similar so we just consider the space $(\mathcal{H}(c_{0}),\tau_{\omega})$.

 In view of (9) it suffices, by Theorem 1,  to take an arbitrary bounded polydisc
 $$A:=\{(z_{m})_{m=1}^{\infty}:|z_{m}|\leq \lambda_{m} {\textrm{ all }} m\}$$ and to show that
 the basis constant, $c_{n}$,
 for $(\mathcal{P}(^{n}c_{0}),\|\cdot\|_{A})$ satisfies $c_{n}\leq 3^{n}$ for all $n$.
 The square ordering on
 $\mathcal{P}^{1}(c_{0})=c_{0}^{\prime}=\ell_{1}$ is just the
 standard ordering and the basis constant is 1. We now suppose
 that $c_{n}\leq 3^{n}$ and prove the required result by induction.

 Let $(\alpha_{m})_{m=1}^{\infty}$ denote an arbitrary sequence of
 scalars. Fix positive integers $s$ and $t$, $s<t$.
 For some non-negative integer $k$ we have the expansion
\begin{eqnarray}\sum_{m=1}^{s}\alpha_{m}P_{n+1,m}&=&\sum_{u=1}^{k+1}
\big{\{}\sum_{m,l(P_{n+1,m})=u}
\alpha_{m}P_{n+1,m}\big{\}}.\end{eqnarray} Letting $Q_{n+1,u}=
\sum_{1\leq m \leq s,l(P_{n+1,m})=u} \alpha_{m}P_{n+1,m}$ for
$1\leq u\leq k+1$ we obtain
\begin{eqnarray*}\sum_{m=1}^{s}\alpha_{m}P_{n+1,m}&=&\sum_{u=1}^{k+1}Q_{n+1,u}\end{eqnarray*}
Note that each $P_{n+1,m}$ is a monomial of degree $n+1$ and that
if $m_{1}<m_{2}$ then $$l(P_{n+1,m_{1}})\leq l(P_{n+1,m_{2}}).$$
If $s<m\leq t$ then $l(P_{n+1,m})\geq k+1$ and for some integer
$k^{*}\geq k+1$
\begin{eqnarray*}\sum_{m=s+1}^{t}\alpha_{m}P_{n+1,m}=\sum_{u=k+1}^{k^{*}}
\big{\{}\sum_{s<m\leq t,l(P_{n+1,m})=u}
\alpha_{m}P_{n+1,m}\big{\}}.\end{eqnarray*} If
$$Q_{n+1,u}^{*}= \sum_{s<m\leq t,l(P_{n+1,m})=u}
\alpha_{m}P_{n+1,m}$$ for $k+1< u\leq k^{*}$, then with the
convention $\sum_{u>k+1}^{k^{*}}=0$ when $k+1=k^{*}$, we have

\begin{eqnarray}\sum_{m=s+1}^{t}\alpha_{m}P_{n+1,m}&=&\sum_{m>s,l(P_{n+1,m})=k+1}
\alpha_{m}P_{n+1,m}+\sum_{u>k+1}^{k^{*}}Q_{n+1,u}^{*}.\end{eqnarray}
If we let
\begin{eqnarray}Q_{n+1,k+1}^{*}=Q_{n+1,k+1}+\sum_{s<m\leq
t,l(P_{n+1,m})=k+1} \alpha_{m}P_{n+1,m}\end{eqnarray} then
\begin{eqnarray*}\sum_{m=1}^{t}\alpha_{m}P_{n+1,m}
&=&\sum_{m=1}^{s}\alpha_{m}P_{n+1,m}+\sum_{m=s+1}^{t}\alpha_{m}P_{n+1,m}\\
&=& \sum_{u=1}^{k+1}Q_{n+1,u}+\sum_{s<m\leq t,l(P_{n+1,m})=k+1}
\alpha_{m}P_{n+1,m}+\sum_{u>k+1}^{k^{*}}Q_{n+1,u}^{*}
\\&=& \sum_{u=1}^{k}Q_{n+1,u}+Q_{n+1,k+1}^{*}+\sum_{u>k+1}^{k^{*}}Q_{n+1,u}^{*}
\\&=& \sum_{u=1}^{k}Q_{n+1,u}+\sum_{u=k+1}^{k^{*}}Q_{n+1,u}^{*}.\end{eqnarray*}
 This identity and (7) imply
 \begin{eqnarray}\|\sum_{u=1}^{k}Q_{n+1,u}\|_{A}&\leq &\|\sum_{m=1}^{t}\alpha_{m}P_{n+1,m}\|_{A}.
 \end{eqnarray} If $l(P_{n+1,m})\leq k$ for $m<m_{0}$ and $l(P_{n+1,m_{0}})=k+1$
 then, by (10) and (12), \begin{eqnarray}Q_{n+1,k+1}=\sum_{m=m_{0}}^{s}\alpha_{m}P_{n+1,m}
 =e_{k+1}^{*}\cdot \sum_{m=m_{0}}^{s}\alpha_{m}P_{n,m-m_{0}+1}\end{eqnarray}
and
\begin{eqnarray}Q_{n+1,k+1}^{*}=\sum_{m=m_{0}}^{s^{*}}\alpha_{m}P_{n+1,m}
 =e_{k+1}^{*}\cdot\sum_{m=m_{0}}^{s^{*}}\alpha_{m}P_{n,m-m_{0}+1}\end{eqnarray}
 for some integer $s^{*}$, $s\leq s^{*}\leq t$. Applying in turn (14), induction, (15),
 and (13), we obtain
\begin{eqnarray*}\|Q_{n+1,k+1}\|_{A}&=&\|e_{k+1}^{*}\cdot
\sum_{m=m_{0}}^{s}\alpha_{m}P_{n,m-m_{0}+1}\|_{A}\\&=&\|e_{k+1}^{*}\|_{A}\cdot
\|\sum_{m=m_{0}}^{s}\alpha_{m}P_{n,m-m_{0}+1}\|_{A}\\&\leq
&3^{n}\|e_{k+1}^{*}\|_{A}\cdot
\|\sum_{m=m_{0}}^{s^{*}}\alpha_{m}P_{n,m-m_{0}+1}\|_{A}\\&=&3^{n}\|e_{k+1}^{*}\cdot
\sum_{m=m_{0}}^{s^{*}}\alpha_{m}P_{n,m-m_{0}+1}\|_{A}\\&=
&3^{n}\|\sum_{m=m_{0}}^{s^{*}}\alpha_{m}P_{n+1,m}\|_{A}\\&=
&3^{n}\|Q_{n+1,k+1}^{*}\|_{A}\\&\leq
&3^{n}\|\sum_{u=k+1}^{k^{*}}Q_{n+1,u}^{*}\|_{A}\\&\leq
&3^{n}\big{\{}\|\sum_{u=1}^{k}Q_{n+1,u}+\sum_{u=k+1}^{k^{*}}Q_{n+1,u}^{*}\|_{A}+
\|\sum_{u=1}^{k}Q_{n+1,u}\|_{A}\big{\}}\\&\leq &2\cdot
3^{n}\|\sum_{m=1}^{t}\alpha_{m}P_{n+1,m}\|_{A}.\end{eqnarray*}
This estimate, together with (7) and (13), implies
\begin{eqnarray*}\|\sum_{m=1}^{s}\alpha_{m}P_{n+1,m}\|_{A}
&=&\|\sum_{u=1}^{k+1}Q_{n+1,u}\|_{A}\\&\leq
&\|\sum_{u=1}^{k}Q_{n+1,u}\|_{A}+\|Q_{n+1,k+1}\|_{A}\\&\leq
&(1+2\cdot 3^{n})\|
 \sum_{m=m_{0}}^{t}\alpha_{m}P_{n+1,m}\|_{A}\end{eqnarray*} and hence
 $$c_{n+1}\leq 1+2\cdot 3^{n}\leq 3^{n+1}.$$  This completes the proof.

sean.dineen@ucd.ie\\School of Mathematical Sciences,\\University
College Dublin,\\Dublin 4,
Ireland.\\\\mujica@ime.unicamp.br\\IMECC-UNICAMP,
\\Rua Sergio Buarque de Holanda 651,\\13083-859 Campinas,\\ SP,
Brazil.
 \end{document}